\documentclass{amsart}
\usepackage{amsmath, amssymb}
\newtheorem{Lem}{Lemma}
\newtheorem{Them}{Theorem}
\newtheorem{Con}{Conjecture}

\newcommand{\fq}{\mathbb{F}_q}
\newcommand{\rmv}[1]{}

\begin{document}
\title{On a conjecture of polynomials with prescribed range}
\author[Muratovi\'{c}-Ribi\'{c}]{Amela Muratovi\'{c}-Ribi\'{c}}
\address{University of Sarajevo, Department of Mathematics, Zmaja od Bosne 33-35, 71000 Sarajevo, Bosnia and Herzegovina}
\email{amela@pmf.unsa.ba}

\author[Wang]{Qiang Wang}
\thanks{Research is partially supported by NSERC of Canada.}
\address{School of Mathematics and Statistics,
Carleton University, Ottawa, Ontario\\ $K1S$ $5B6$, CANADA}

\email{wang@math.carleton.ca}

\begin{abstract}
We show that, for any integer $\ell$ with $q-\sqrt{p} -1 \leq \ell <
q-3$ where $q=p^n$ and $p>9$, there exists a multiset $M$ satisfying
that $0\in M$ has the highest multiplicity  $\ell$ and $\sum_{b\in
M} b =0$ such that every polynomial over finite fields $\fq$ with
the prescribed range $M$ has degree greater than $\ell$. This
implies that Conjecture 5.1. in \cite{gac} is false over finite
field $\fq$ for  $p > 9$  and $k:=q-\ell -1 \geq 3$.
\end{abstract}

\maketitle
\section{Introduction}

Let $\fq$ be a finite field of $q = p^n$ elements and $\fq^*$ be the
set of all nonzero elements. Any mapping from $\fq$ to itself can be uniquely
represented by a polynomial of degree at most $q-1$. The degree of
such a polynomial is called the {\it reduced} degree.  A multiset $M$ of size $q$ of field elements is called 
the {\it range} of the  polynomial $f(x) \in \fq[x]$ if $M = \{ f(x) : x \in \fq \}$ as a multiset
(that is, not only values, but also multiplicities need to be the
same). Here we use the set notation for multisets as well. We refer the readers to \cite{gac} for more details.  In the study of polynomials
with prescribed range, G\`{a}cs et al. recently proposed the
following conjecture.
\begin{Con} [Conjecture 5.1, \cite{gac}]\label{conjecture1}
Suppose $M=\{a_1,a_2,\ldots ,a_q\}$ is a multiset of $\fq$
 with $a_1+\ldots +a_q=0$, where $q=p^n$, $p$ prime. Let
$k<\sqrt{p}$. If there is no polynomial with range $M$ of degree
less than $q-k$, then $M$ contains an element of multiplicity at
least $q-k$.
\end{Con}

We note that Conjecture~\ref{conjecture1}  is equivalent to
\begin{Con}\label{contrapositive}
 Suppose $M=\{a_1,a_2,\ldots ,a_q\}$ is a multiset of $\fq$
 with $a_1+\ldots +a_q=0$, where $q=p^n$, $p$ prime. Let
$k<\sqrt{p}$. If multiplicities of all elements in $M$ are less than
$q-k$, then there exist a polynomial with range $M$ of the degree
less than $q-k$.
\end{Con}

In the case $k=2$,  Conjecture~\ref{conjecture1} holds by
Theorem~2.2.\ in \cite{gac}.  In particular, Theorem~2.2 in
\cite{gac} gives a complete description of $M$ so that there is no
polynomial with range $M$ of reduced degree less than $q-2$.  In
this paper, we study the above conjecture for $k \geq 3$.

Suppose we take a prescribed range $M$ such that the highest
multiplicity in $M$ is $\ell =q-k-1$, if the above conjecture were
true then it follows that there exist a polynomial, say $g(x)$, with
range $M$  and the degree of $g(x)$ is less than $q-k$. On the other
hand,   If $a\in M$ is the element with multiplicity $\ell$ then
polynomial $g(x)-a$ has $\ell$ roots and thus
 the degree of $g(x)$ is at least equal to the
highest multiplicity $\ell$ in $M$.  Therefore the degree of $g(x)$
must be $\ell =q-k-1$.  This means that,  if
Conjecture~\ref{contrapositive} were true, then for every multiset
$M$ with the highest multiplicity $\ell =q-k-1$ where $1\leq k<
\sqrt{p}$ there exists a polynomial with range $M$ of the degree
$\ell $. Note that $k<\sqrt{p}$  implies $\ell
=q-k-1>q-\sqrt{p}-1\geq \frac{q}{2}$ when $q > 5$. Also $3\leq k\leq
\sqrt{p}$ implies $p>9$.

Let $M=\{a_1,a_2,\ldots a_q\}$ be  a given multiset. We  consider
polynomials $f(x):\fq \rightarrow M$, with the least degree. Denote
by $\ell $ the highest multiplicity in $M$ and let $\ell +m=q$.  If
$a\in M$ is an element with multiplicity $\ell $ then the polynomial
$f(x)-a$ has the same degree as $f(x)$ and $0$ is in the range of
$f(x)-a$ such that $0$ has the same highest multiplicity $\ell$.
Therefore, we only consider multisets $M$ where $0$ has the highest
multiplicity for the rest of paper.

In particular, we prove  the following theorem.

\begin{Them} \label{main}
Let $\fq$ be a finite field of $q =p^n$ elements with $p > 9$. For
every $\ell$ with $q-\sqrt{p}-1\leq \ell < q-3$ there exists a
mutiset $M$ with $\sum_{b\in M} b =0$ and the highest mutiplicity
$\ell$ achieved at $0\in M$ such that every polynomial over the
finite field $\fq$ with the prescribed range $M$ has degree greater
than $\ell $.
\end{Them}


In particular,  for any $p > 9$, if  we take $\ell =q-k-1\leq q-4$,
i.e., $k\geq 3$, then Theorem~\ref{main} implies that
Conjecture~\ref{contrapositive} fails.

\section{Proof of Theorem~\ref{main}}
In this section, we prove Theorem~\ref{main}. Let $\ell $ be fixed
and $q-\sqrt{p}-1 \leq \ell <q-3$. Because $p > 9$,  $\sqrt{p} > 3$
and such $\ell$ exists.
 Let $M$ be a multiset such that
$0\in M$  has the highest multiplicity $\ell $ and $\sum_{b\in
M}b=0$. Note that $\ell \geq \frac{q}{2}$ implies that multipilicity
of any nonzero element in $M$ is less than $m := q-\ell \leq
\frac{q}{2}$ (indeed the highest multiplicity is achieved at $0$).
Let $f:\fq \rightarrow M$. Let $U\subseteq \fq$ such that
$f(U)=\{0^\ell \}$ (the multiset of $\ell$ zeros) and $T=\fq \setminus U$,  i.e.,  $x\in T$ implies
$f(x)\neq 0$. Then $|U|=\ell $ and $|T|=m$ and $M=f(U)\cup f(T)$.
Then polynomial $f:\fq\rightarrow M$ can be written in the form
$f(x)=h(x)P(x)$ where $P(x)=\prod_{s\in U}(x-s)$ and  $h(x)\neq 0$
has no zeros in $T$. Then $\deg(f)\geq \deg(P) =  \ell $. We note
that there is a bijection between polynomials with range $M =\{a_1, \ldots, a_q\}$ and the
ordered sets $(b_1, \ldots, b_q)$ (that is, permutations) of $\fq$: a
permutation corresponds  to the function $f(b_i) = a_i$. For each
$U$, there are many different $h(x)$'s corresponding to different
ordered sets $(b_1, \ldots, b_q)$ such that $f(b_i) =0$ for all $b_i
\in U$.  However, if  $h(x)=\lambda\in \fq^*$ then $f(x)$ is a
polynomial of the least degree  and each polynomial $f(x)$ is
uniquely determined by a set $T$ and a nonzero scalar $\lambda $.
Thus we denote $f(x)$ by
\begin{equation} \label{poly}
f_{(\lambda ,T)}(x)=\lambda \prod_{s\in \fq\setminus T}(x-s).
\end{equation}

Therefore its range $M$ is also uniquely determined by $T$ and
$\lambda $. Denote by $\mathcal{T}$ the family of all subsets of
$\fq $ of cardinality $m$, i.e.,
$$\mathcal{T}=\{T \mid T\subseteq \fq, |T|=m\}.$$
Denote by $\mathcal{M}$ the family of all multisets $M$ of order $q$
containing $0$, having the highest multiplicity $\ell $ achieved at
$0$ and whose sum of elements in $M$ is equal to the $0$, i.e.,
$$
\mathcal{M}=\{M \mid 0\in M, \text{ multiplicity}(0)=\ell,
\sum_{b\in M}b=0\}.
$$

Equation~(\ref{poly}) uniquely determines a mapping
$$\mathcal{F}:\fq^* \times \mathcal{T}\rightarrow \mathcal{M}$$
where $$(\lambda ,T)\mapsto \text{range}(f_{\lambda, T}(x)).$$

 Now  by Equation~(\ref{poly}) it follows
that for every $\hat{s}\in T$ we have
\begin{equation}\label{basicequations}
f_{\lambda, T}(\hat{s})=\lambda P(\hat{s})=\lambda \prod_{s\in
\fq,s\neq \hat{s}}(\hat{s}-s)\Bigl(\prod_{s\in T,s\neq
\hat{s}}(\hat{s}-s)\Bigr)^{-1}=-\lambda \Bigl(\prod_{s\in T,s\neq
\hat{s}}(\hat{s}-s)\Bigr)^{-1}.\end{equation} (Note that this
equation does not hold for $x\in \fq\setminus T$). In the following
we find an upper bound of $| \text{range}(\mathcal{F}) |$ and a
lower bound of $| \mathcal{M} |$ and show that $|\mathcal{M}| > |\text{range}(\mathcal{F})|$. This implies that
Theorem~\ref{main} holds.

First of all we observe

\begin{Lem} Let $\lambda $ and $T$ be given. For any $c\in
\fq^*$ and any $b\in \fq$, we have
$$f_{(\lambda,T)}(\hat{s})=f_{(c^{m-1}\lambda, cT+b)}(c\hat{s}+b),\qquad {\text for}~\hat{s}\in T$$
i.e.,
$$\mathcal{F}(\lambda, T)=\mathcal{F}(c^{m-1}\lambda ,cT+b).$$

\end{Lem}
\begin{proof}
We use notation $cT+b=\{cs+b \mid s\in T\}$. Substituting in
(\ref{basicequations}), we obtain $f_{(c^{m-1}\lambda,
cT+b)}(c\hat{s}+b)=-c^{m-1}\lambda (\prod_{s\in T,s\neq
\hat{s}}((c\hat{s}+b)-(cs+b)))^{-1}=-\lambda (\prod_{s\in T,s\neq
\hat{s}}(\hat{s}-s))^{-1}=f_{(\lambda, T)}(\hat{s})$.
\end{proof}

Now we  use Burnside's Lemma to find an upper bound of the
cardinality of $\text{range}(\mathcal{F})$.

\begin{Lem}\label{countingRange} Let $m< \sqrt{p}+1$ and let $d = \gcd(q-1, m-1)$. Then
$$|\text{range}(\mathcal{F})|\leq \frac{(q-1)(q-2)\ldots
(q-m+1)}{m!}+\sum_{{i\mid d \atop i >1}
}\phi(i)\binom{\frac{q-1}{i}}{\frac{m-1}{i}}.$$
\end{Lem}
\begin{proof}
Let $\mathcal{G}$ be group of all (nonzero) linear polynomials in
$\fq[x]$ with the composition operation. Indeed, $\mathcal{G}$ is a
subgroup of the group of all permutation polynomials because the
composition of two linear polynomials is again a linear polynomial,
the identity mapping is a linear polynomial, and the inverse of a
linear polynomial is again a linear polynomial. We  use notation
$cT+b=\{cs+b \mid s\in T\}$ again. Then $\mathcal{G}$  acts on the
set $\fq^*\times \mathcal{T}$ with $\Phi :
\mathcal{G}\times(\fq^*\times \mathcal{T}) \rightarrow \fq^*\times
\mathcal{T}$, where
$$\Phi : (cx+b, (\lambda,T))\mapsto ( c^{m-1}\lambda ,cT+b).$$
The elements of the same orbit
$$\mathcal{G}(\lambda, T)=\{(c^{m-1}\lambda ,cT+b) \mid cx+b\in\mathcal{G}\}$$
are all mapped to the same element $M\in \mathcal{M}$ by Lemma 1. By
Burnside's Lemma the number of orbits $N$ is given by
$$ N =\frac{1}{|\mathcal{G}|}\sum_{g\in \mathcal{G}}|(\fq^*\times
\mathcal{T})_g|,$$ where  $g(x)=cx+b$, and

$$(\fq^*\times \mathcal{T})_g=\{(\lambda ,T) \mid (\lambda ,T)\in
\fq^*\times \mathcal{T}, (c^{m-1}\lambda ,cT+b)=(\lambda ,T)\}.$$

The equation $cx+b=x$ over $\fq$ is equivalent to $(c-1)x=-b$, which
has exactly one solution if $c\neq 1$; no solutions if $c=1$ and
$b\neq 0$; $q$ solutions if $c=1$ and $b=0$. If $c\neq 1$ and $i
:=ord(c)\mid q-1$, then  this linear polynomial has one fixed
element and  $\frac{q-1}{i}$ cycles of length $i$. If $c=1$ and
$b\neq 0$ then $g^p(x)=x+pb=x$ and thus $g(x)$ has cycles of length
$p$ where $p = char(\fq)$.

Assume $T=cT+b$.  Let $s\in T$. Then $g(s)\in cT+b=T$. So the  cycle
$s,g(s),g^2(s),\ldots ,g^i(s)=s$ is contained in $T$.

This means that, under the assumptions of  $c\neq 1$ and $T=cT+b$,
either $T$ has one fixed element and $\frac{m-1}{i}$ cycles of the
length $i$ which are defined by permutation $g(x)$, or $T$ has
$\frac{m}{i}$ cycles length $i$ which are defined by permutation
$g(x)$. In the latter case, the fixed element of $g(x)$ is in $\fq
\setminus T$.

In the former case, if  $c\in \fq^* \setminus \{1\}$ satisfies $i=ord(c) \mid
d = \gcd(q-1,m-1)$ then there are
$\binom{\frac{q-1}{i}}{\frac{m-1}{i}}$ sets fixed by $g(x)$.
Moreover,  $c^{m-1}=(c^i)^\frac{m-1}{i}=1$. Hence, for each set $T$
fixed by $g(x)$ and any $\lambda \in \fq^*$ we must have
$(c^{m-1}\lambda ,cT+b)=(\lambda,T)$. This implies that
$$|(\fq^*\times
\mathcal{T})_g|=(q-1)\binom{\frac{q-1}{i}}{\frac{m-1}{i}}.$$

If $c\in \fq^*$ satisfies $i = ord(c) \mid \gcd(q-1,m)$ then there
are $\binom{\frac{q-1}{i}}{\frac{m}{i}}$ sets $T$ fixed by $g(x)$.
But for each $T$ fixed by $g(x)$, $c^{m-1}=c^{-1}\neq 1$ and thus
$(c^{m-1}\lambda ,cT+b)\neq (\lambda, T)$. Therefore
$$|(\fq^*\times
\mathcal{T})_g|=0.$$

If $c=1$ and $b=0$ then $g(x) =x$. So $|(\fq^*\times
\mathcal{T})_g|=(q-1)\binom{q}{m}$.  If $c=1$ and $b\neq 0$ then
$cT+b \neq T$. Otherwise, it  implies that $T$ contains elements of
the cycles of the length $p$ which  contradicts to $m<\sqrt{p}+1$.

Since $d = \gcd(q-1, m-1)$,  we obtain

\begin{eqnarray*}
N&=&\frac{1}{|\mathcal{G}|}\sum_{g\in \mathcal{G}}|(\fq^*\times
\mathcal{T})_g|\\
&=&\frac{1}{q(q-1)}\biggl( (q-1)\cdot
\binom{q}{m}+\sum_{\begin{array}{c}\scriptscriptstyle c\in \fq^*\setminus \{1\}\\
\scriptscriptstyle i=ord(c)\mid
d\\
\scriptscriptstyle b\in \fq \end{array}} (q-1)\binom{\frac{q-1}{i}}{\frac{m-1}{i}} \biggr),\\
&=&\frac{1}{q(q-1)}\biggl((q-1)\cdot
\binom{q}{m}+q(q-1)\sum_{\begin{array}{c} \scriptscriptstyle c\in \fq^*\setminus \{1\}\\
\scriptscriptstyle i=ord(c)\mid  d
\end{array}} \binom{\frac{q-1}{i}}{\frac{m-1}{i}} \biggr)\\
&=&\frac{(q-1)(q-2)\ldots (q-m+1)}{m!}+\sum_{i\mid
d}\phi(i)\binom{\frac{q-1}{i}}{\frac{m-1}{i}},
\end{eqnarray*}
where  $\phi(i)$ is the number of $c$'s such that the order of $c$
is $i > 1$.

 Since two orbits could  possibly be mapped to the same
multiset $M\in \mathcal{M}$ we finally have an inequality
\begin{equation}|\text{range}(\mathcal{F})|\leq \frac{(q-1)(q-2)\ldots
(q-m+1)}{m!}+\sum_{i\mid
d}\phi(i)\binom{\frac{q-1}{i}}{\frac{m-1}{i}}.\end{equation}
\end{proof}

\par Now we find a lower bound of the cardinality of
$\mathcal{M}=\{\stackrel{l\quad \text{ times}}{\overbrace{0,0\ldots
,0}},b_1,b_2,\ldots ,b_m\}$ such that  $b_i \neq 0$ for $i=1,
\ldots, m$ and
\begin{equation}\label{imagesum}
b_1+b_2+\ldots +b_m=0.
\end{equation}
Although we can find a simpler exact formula for the number of solutions to Equation~(\ref{imagesum}), we prefer the following lower bound for $|\mathcal{M}|$ which has the same format as the upper bound  of $|\text{range}(\mathcal{F})|$ in order to compare them directly. 

\begin{Lem}\label{countingM}
Let $A=1$ if $m-1\mid q-1$ and $A=0$ otherwise. If $m\geq 6$ then
$$
|\mathcal{M}|\geq \frac{(q-1)\ldots
 (q-m+2)(q-2)}{m!}+$$
 $$\sum_{\begin{matrix}
\scriptscriptstyle 1<i<m-1\\
\scriptscriptstyle i\mid
\gcd(q-1,m-1)\end{matrix}}\frac{[(q-1)\ldots
(q-\frac{m-1}{i}+2)][(q-\frac{m-1}{i}+1)(q-\frac{m-1}{i}-1)+\frac{m-1}{i}]}{\Bigl(\frac{m-1}{i}\Bigr)!}+A(q-1).$$
  If $m=4$  and $3\mid q-1$ then
$$|\mathcal{M}|\geq \frac{(q-1)(q-2)^2}{4!}.$$
If $m=5$ then
$$|\mathcal{M}|\geq \frac{(q-1)(q-2)^2(q-3)}{5!}+ A(q-1).$$
\end{Lem}
\begin{proof}
In order to give a lower bound of $|\mathcal{M}|$, we count two
different classes of  families of multisets $M$. The first class
contains families of those multisets $M$ such that almost all
nonzero elements $b_i$'s have the same multiplicities  greater than
one except the last element $b_m$. And the second family  class
contains those multisets $M$ such that almost all nonzero elements
$b_i$'s have multiplicities one except that the last two elements
$b_{m-1}$ and $b_{m}$.

First,  we count those multisets $M$ such that almost all nonzero
elements $b_i$'s have the same multiplicities  greater than one
except the last element $b_m$.  That is, for any $i$ such that
$1<i<m-1$ and $i\mid \gcd(q-1,m-1)$,
 we want to choose $\frac{m-1}{i}$ pairwise distinct nonzero
elements each of multiplicity $i$ so that
$\displaystyle\sum_{j=1}^{\frac{m-1}{i}}ib_j\neq 0$ (the sum being
equal to zero would imply $b_m=0$, a contradiction).  For each such
$i$, we denote the family of these multisets by $\mathcal{M}_i$.

We note that each multiset $M\in \mathcal{M}_i$  can be written as
$$M=\{\stackrel{\ell \quad \text{ times}}{\overbrace{0,0,\ldots
,0}},\stackrel{i\quad \text{ times}}{\overbrace{b_1,\ldots
,b_1}},\ldots ,\stackrel{i\quad \text{
times}}{\overbrace{b_{\frac{m-1}{i}},\ldots
,b_{\frac{m-1}{i}}}},b_m\}.$$

Obviously each  multiset is invariant to the ordering.  However, let
us first consider  the ordered tuples $(b_1,\ldots
,b_{\frac{m-1}{i}-1})$ satisfying that $b_i$'s are nonzero and
pairwise distinct.  Out of a total of $(q-1)\ldots
(q-\frac{m-1}{i}+1)$ such ordered tuples, there are
 $(q-1)\ldots
(q-\frac{m-1}{i}+2)\frac{m-1}{i}$ choices such that
$-\displaystyle\sum_{j=1}^{\frac{m-1}{i}-1}b_j\in \{0,b_1,\ldots
,b_{\frac{m-1}{i}-1} \}$ and $(q-1)\ldots (q-\frac{m-1}{i}+2)(q- 2
\frac{m-1}{i}+1)$ ordered tuples such that
$-\displaystyle\sum_{j=1}^{\frac{m-1}{i}-1}b_j\not \in
\{0,b_1,\ldots ,b_{\frac{m-1}{i}-1} \}$.  If
$-\displaystyle\sum_{j=1}^{\frac{m-1}{i}-1}b_j\in \{0,b_1,\ldots
,b_{\frac{m-1}{i}-1} \}$  then $b_{\frac{m-1}{i}}$ can be chosen in
$ q-\frac{m-1}{i}$ ways and otherwise it can be chosen in
$q-\frac{m-1}{i}-1$ way. Because the element $b_m$ is uniquely
determined by $i\sum_{j=1}^{\frac{m-1}{i}}b_j$, we have in total
$$(q-1)\ldots
(q-\frac{m-1}{i}+2)\frac{m-1}{i}(q-\frac{m-1}{i})$$
$$+ (q-1)\ldots
(q-\frac{m-1}{i}+2)(q-2\frac{m-1}{i}+1)(q-1-\frac{m-1}{i})$$
$$= \left( (q-1)\ldots(q-\frac{m-1}{i}+2)\right) \left((q-\frac{m-1}{i}+1)(q-\frac{m-1}{i}-1)+\frac{m-1}{i}\right)$$
ordered tuple $(b_1\ldots ,b_{\frac{m-1}{i}})$ satisfying
Equation~(\ref{imagesum}) and that $b_i \neq 0$ for $i=1, \ldots, m$
and each element is of multiplicity $i$ except that last element.

Since there are $\Bigl(\frac{m-1}{i}\Bigr)!$ permutations of the
ordered tuples $(b_1,\ldots ,b_{\frac{m-1}{i}})$, there are
$$\frac{[(q-1)\ldots(q-\frac{m-1}{i}+2)][(q-\frac{m-1}{i}+1)(q-\frac{m-1}{i}-1)+\frac{m-1}{i}]}{\Bigl(\frac{m-1}{i}\Bigr)!}$$
 elements in
$\mathcal{M}_i$.

Similarly, if $m-1\mid q-1$, we denote by $\mathcal{M}_{m-1}$ the
set of multisets $M$ such that all $b_i$'s are the same nonzero
element for $i=1, \ldots, m-1$ and their sum together with $b_m$ is
zero. It is easy to see that there are $q-1$ such $M$'s,  i.e.,
$|\mathcal{M}_{m-1}|=q-1$.

Now we show that $\mathcal{M}_i\cap \mathcal{M}_j =\emptyset$ for $1
< i\neq j \leq m-1$. We prove this by contradiction and  we use
heavily the fact that, for each $i$, there are $\frac{m-1}{i}+1$
distinct elements in $M\in \mathcal{M}_i$ if $b_m\neq b_k$ for
$1\leq k\leq \frac{m-1}{i}$ and
 there are $\frac{m-1}{i}$ distinct elements in $M$ if   $b_m = b_k$ for some $k$.  Assume that $\mathcal{M}_i\cap \mathcal{M}_j \neq \emptyset$. Obviously,
  $\frac{m-1}{i}\neq \frac{m-1}{j}$ because $i \neq j$. Hence
either $\frac{m-1}{i}+1 = \frac{m-1}{j}$ or $\frac{m-1}{j}+1 =
\frac{m-1}{i}$.

Without loss of generality, we assume now $\frac{m-1}{i}+1=
\frac{m-1}{j}$ and $M\in \mathcal{M}_i\cap \mathcal{M}_j$. Then in
the multiset $M$, we have $\frac{m-1}{i}$ elements of multiplicity
$i$ and one element of multiplicity $1$ since $M\in \mathcal{M}_i$.
Moreover,  the number of elements of multiplicity $j$ is
$\frac{m-1}{j}-1$ and there is one element of multiplicity $j+1$
since $M\in \mathcal{M}_j$. Because $i >j$, we must have $i=j+1$ and
$j=1$ by comparing the multipicities. However,  this implies we must
have $\frac{m-1}{i}=1$ and $\frac{m-1}{j}-1=1$. Hence $i = m-1$ and
$j=\frac{m-1}{2}$, contradicts to $i =j+1$ when $m > 3$.

Therefore $\mathcal{M}_i\cap \mathcal{M}_j\neq \emptyset $ for all $
1 < i\neq j \leq m-1$. Now for $m \geq 4$ we have
$$|\bigcup_{\begin{matrix}\scriptscriptstyle 1<i\leq m-1\\
\scriptscriptstyle i\mid \gcd(q-1,m-1)\end{matrix}}\mathcal{M}_i|=
A(q-1) +$$
$$\sum_{\begin{matrix}\scriptscriptstyle 1<i<m-1\\
\scriptscriptstyle i\mid
\gcd(q-1,m-1)\end{matrix}}\frac{[(q-1)\ldots
(q-\frac{m-1}{i}+2)][(q-\frac{m-1}{i}+1)(q-\frac{m-1}{i}-1)+\frac{m-1}{i}]}{\Bigl(\frac{m-1}{i}\Bigr)!}.$$

Next we count those multisets $M$ such that almost all nonzero
elements $b_i$'s have multiplicities one except that the last two
elements $b_{m-1}, b_{m}$. That is, $b_1,\ldots, b_{m-2}$ are
pairwise distinct nonzero elements, $b_{m-1} \neq 0$ is chosen in a
way such that $\displaystyle\sum_{j=1}^{m-1}b_j\neq 0$,
 and $b_m$ is uniquely
determined by $\displaystyle\sum_{j=1}^mb_j=0$. The family of such
multisets is denoted by $\mathcal{M}_0$. We note that $b_{m-1}$ and
$b_m$ could be same as one of $b_j$'s where $j=1, \ldots, m-2$. So
the highest mulitiplicity is at most $3$.

Consider all $(q-1)\ldots (q-m+2)$ different ordered tuples
$(b_1,\ldots ,b_{m-2})$.  If $-\displaystyle\sum_{j=1}^{m-2}b_j\neq
0$
 we can choose $b_{m-1}$ in $q-2$ ways and
otherwise there are $q-1$ choices for $b_{m-1}$. Thus in total there
are at least $(q-1)\ldots (q-m+2)(q-2)$
 ordered tuples $(b_1,\ldots ,b_m)$.

Let $S_1$ be the  number of such ordered tuples without repetition,
$S_2$ be the number of ordered tuples with exactly one repeated
element, $S_3$ be the number of arrays with exactly two pairs of
repeated elements, and $S_4$ be the number of tuples with exactly
one element repeated $3$ times. Because multisets are invariant to
the ordering, there are  at least
$$\frac{S_1}{m!}+\frac{S_2}{(m-1)!}+\frac{S_3}{(m-2)!2!}+\frac{S_4}{(m-2)!}\geq\frac{(q-1)\ldots
 (q-m+2)(q-2)}{m!}$$
 such multisets in $\mathcal{M}_0$, i.e.,  $$|\mathcal{M}_0|\geq \frac{(q-1)\ldots
 (q-m+2)(q-2)}{m!}.$$

 We note that each multiset from $\mathcal{M}_0$ contains at least $m-2$ distinct  elements and  each multiset from $\mathcal{M}_i$ with $i>1$ contains at most $\frac{m-1}{i}+1\leq \frac{m-1}{2}+1$ distinct
 elements. Since $\frac{m-1}{2}+1<m-2$ for $m\geq 6$ we have that
 $\mathcal{M}_0\cap\mathcal{M}_i=\emptyset $ as long as $m\geq 6$.
Therefore we can conclude that for $m\geq 6$ we have
$$|\mathcal{M}|\geq |\mathcal{M}_0|+|\bigcup_{\begin{matrix}\scriptscriptstyle 1<i\leq m-1\\
\scriptscriptstyle i\mid
\gcd(q-1,m-1)\end{matrix}}\mathcal{M}_i|\geq \frac{(q-1)\ldots
 (q-m+2)(q-2)}{m!}+$$
 $$\sum_{\begin{matrix}\scriptscriptstyle 1<i<m-1\\
\scriptscriptstyle i \mid
\gcd(q-1,m-1)\end{matrix}}\frac{[(q-1)\ldots
(q-\frac{m-1}{i}+2)][(q-\frac{m-1}{i}+1)(q-\frac{m-1}{i}-1)+\frac{m-1}{i}]}{\Bigl(\frac{m-1}{i}\Bigr)!}+A(q-1).$$

Let $m=4$. If $i>1$ and $i\mid
 \gcd(m-1,q-1)$ then $i=3$. Thus in this case
 $\mathcal{M}_3\cap \mathcal{M}_0=\{ \{a,a,a,b\} \mid a\in \fq^*, b=-3a\neq
 a \}$ since $p> 9$.
 By the principle of the inclusion-exclusion  we obtain
$$|\mathcal{M}|\geq \frac{(q-1)(q-2)(q-2)}{4!}+(q-1)-(q-1)=\frac{(q-1)(q-2)(q-2)}{4!}.$$

If $m=5$, then $i>1$ and $i\mid \gcd(4, q-1)$ imply $i=2$ or $i=4$.
Obviously $\mathcal{M}_0 \cap \mathcal{M}_4 = \emptyset$ because
each element in a multiset of $\mathcal{M}_0$ has multipliciy at
most $3$.  Similarly, any multiset in both $\mathcal{M}_0$ and
$\mathcal{M}_2$ must contain $\frac{m-1}{i}+1=m-2=3$ distinct
elements, two of them come in pairs. That is,

$$\mathcal{M}_0\cap \mathcal{M}_2=\Bigl\{\{a,a,b,b,c\} \mid a,b,c\in
\fq^*,a\neq b,a\neq c,b\neq c\Bigr\}.$$ If $a$ is chosen in $q-1$
ways  then $b\not\in\{0,a,-a\}$ and we can choose $b$ in $q-3$ ways.
Since multisets are invariant to the ordering we have
$$|\mathcal{M}_0\cap \mathcal{M}_2|=\frac{(q-1)(q-3)}{2!}.$$
Again the principle of inclusion-exclusion implies
\begin{eqnarray*}
|\mathcal{M}| &\geq &   |\mathcal{M}_0|+|\mathcal{M}_2|+|\mathcal{M}_4|- |\mathcal{M}_0 \cap \mathcal{M}_2| \\
&= &  \frac{(q-1)(q-2)(q-3)(q-2)}{5!}+ \frac{(q-1)(q-3)}{2!}+A(q-1)-\frac{(q-1)(q-3)}{2!} \\
&= & \frac{(q-1)(q-2)(q-3)(q-2)}{5!}+A(q-1).
\end{eqnarray*}

\end{proof}

We need the following simple result to compare the bounds of
$\mathcal{M}$ and $|\text{range}(\mathcal{F})|$ in order to complete
the proof of Theorem~\ref{main}.

\begin{Lem} \label{inequality}
(i) For $m\geq 4$,   we have
$$\frac{(q-1)(q-2)\ldots
(q-m+1)}{m!}<\frac{(q-1)\ldots (q-m+2)(q-2)}{m!}.$$

(ii) If $1<i<m-1$ and $i\mid \gcd(q-1,m-1) $ then
$$\phi(i)\binom{\frac{q-1}{i}}{\frac{m-1}{i}}<\frac{(q-1)\ldots
(q-1-\frac{m-1}{i}+2)[(q-\frac{m-1}{i}+1)(q-\frac{m-1}{i}-1)+\frac{m-1}{i}]}{(\frac{m-1}{i})!}.$$

(iii) If $i=m-1\mid q-1$ then
$$\phi(m-1)\binom{\frac{q-1}{i}}{\frac{m-1}{i}}<q-1.$$
\end{Lem}
\begin{proof}
\par (\emph{i}) Clearly, $q-m+1<q-2$ for $m\geq 4$.
\par (\emph{ii}) The inequality
$$\phi(i)\binom{\frac{q-1}{i}}{\frac{m-1}{i}}<\frac{(q-1)\ldots
(q-\frac{m-1}{i}+2)[(q-\frac{m-1}{i}+1)(q-\frac{m-1}{i}-1)+\frac{m-1}{i}]}{(\frac{m-1}{i})!}$$
is equivalent to
$$
\phi(i)\left(\frac{q-1}{i}\bigl(\frac{q-1}{i}-1\bigr)\ldots
\bigl(\frac{q-1}{i}-\frac{m-1}{i}+1\bigr)\right) $$
$$<(q-1)(q-2)\ldots (q-\frac{m-1}{i}+2) \left((q-\frac{m-1}{i}+1)(q-\frac{m-1}{i}-1)+\frac{m-1}{i}\right).$$
Using $\phi(i)\frac{q-1}{i}<q-1,$ $\frac{q-1}{i}-j<\frac{q-1}{i}$
for $j=1,\ldots ,\frac{m-1}{i}-2$ and $\frac{q-1}{i} -\frac{m-1}{i}
+ 1 <  q-1 - \frac{m-1}{i}$ (since $i >1$),  we have
\begin{eqnarray*}
&& \phi(i)\frac{q-1}{i}\bigl(\frac{q-1}{i}-1\bigr)\ldots
\bigl(\frac{q-1}{i}-\frac{m-1}{i}+1\bigr) \\
&<&(q-1)\bigl(\frac{q-1}{i}-1\bigr)\ldots
\bigl(\frac{q-1}{i}-\frac{m-1}{i}+1\bigr) \\
&<& (q-1)(q-2)\ldots
(q-\frac{m-1}{i}+2)(q-\frac{m-1}{i}+1) \bigl(\frac{q-1}{i}-\frac{m-1}{i}+1\bigr)\\
&<& (q-1)(q-2)\ldots
(q-\frac{m-1}{i}+2)(q-\frac{m-1}{i}+1)(q- 1 - \frac{m-1}{i}) \\
&<&(q-1)(q-2)\ldots (q-\frac{m-1}{i}+2) \left(
(q-\frac{m-1}{i}+1)(q-\frac{m-1}{i}-1)+\frac{m-1}{i}\right).
\end{eqnarray*}

\par (\emph{iii}) If $i=m-1\mid q-1$ then $\phi(m-1)\frac{q-1}{m-1}<q-1$.
\end{proof}

\textbf{Proof of Theorem~\ref{main}:}
 If $m\geq 6$ it follows directly from Lemmas~\ref{countingRange}, ~\ref{countingM}, ~\ref{inequality}.

Note that $m\leq \sqrt{p}+1$. If $m=5$ then $5 \leq \sqrt{p} + 1$
implies that $p >16$. Hence we have
\begin{eqnarray*}
|\text{range}(\mathcal{F})| &\leq&
\frac{(q-1)(q-2)(q-3)(q-4)}{5!}+\phi(2)\binom{\frac{q-1}{2}}{2}+A\phi(4)\frac{q-1}{4}\\
& = & \frac{(q-1)(q-2)(q-3)(q-4)}{5!}+\frac{(q-1)(q-3)}{8}+A\frac{q-1}{2}\\
& < & \frac{(q-1)(q-2)(q-3)(q-4)}{5!}+\frac{2(q-2)}{15} \frac{(q-1)(q-3)}{8}+A\frac{q-1}{2}\\
&\leq& \frac{(q-1)(q-2)(q-3)(q-2)}{5!}+A(q-1)\\
&\leq &  |\mathcal{M}|.
\end{eqnarray*}

If $m=4$ and $3\nmid q-1$ then the result follows directly from
Lemmas~\ref{countingRange}, \ref{countingM}, and \ref{inequality}
(i). If $3\mid q-1$ then
$$ |\text{range}(\mathcal{F})|\leq
\frac{((q-1)(q-2)(q-3)}{4!}+\phi(3)\frac{q-1}{3}<\frac{(q-1)(q-2)^2}{4!}<|\mathcal{M}|$$
 holds for $q>18$.  Note that $m\leq \sqrt{p}+1$ implies $p\geq 9$. By the
assumption of $p > 9$,   we must have  $p\geq 11$. 
  The only possible prime power $q \leq 18$ such that $p\geq 11$ and $3 \mid q-1$ is $q=13$. It is easy
  to  compute that the number of all the possible solutions to Equation~(\ref{imagesum}) with desired properites over $\mathbb{F}_{13}$ is $|\mathcal{M}|=105$ by a computer program.
  For $q=13$, then $\gcd(q-1, m-1)= 3$
  and thus
$|\text{range}(\mathcal{F})| \leq 63 < 105 = |\mathcal{M}|$.
 Hence the proof is complete. $\Box$

If $m=2$ and $m=3$ these polynomials satisfying the conjecture do
exist. Indeed, if $m=2$ and
 $b_2=-b_1$, then we can construct the minimum degree polynomial $f(x) = \lambda \prod_{s \in \fq \setminus T} (x-s)$
 with the prescribed range $M=\{0, \ldots, 0, b_1, -b_1\}$  by letting  $T= \{b_1^{-1}, 0\}$ and $\lambda =1$.
\par For the case $m=3$, for any multiset $M = \{0,\ldots ,0,b_1,b_2,b_3\}$ with $b_1+b_2+b_3=0$ such that $b_1, b_2,
b_3$ are all nonzero there exists a polynomial
$f(x)=\lambda\prod_{s\in \fq\setminus T}(x-s)$ of the least degree
with range $M$. Indeed, let $T=\{b_2,-b_1,0\}$ and $\lambda
=b_1b_2b_3$. Then using $b_3=-(b_1+b_2)$ we obtain
$$f(b_2)=b_1b_2b_3\frac{-1}{(b_2+b_1)b_2}=b_1;$$
$$f(-b_1)=b_1b_2b_3\frac{-1}{(-b_1-b_2)(-b_1)}=b_2;$$
$$f(0)=b_1b_2b_3\frac{-1}{(-b_1)(b_2)}=b_3.$$

\vskip 1cm
\end{document}